\documentclass[a4paper,12pt,twoside]{article}
\usepackage[latin1]{inputenc}
\usepackage {lettrine,stmaryrd}
\usepackage[brazil,english]{babel}
\usepackage{amssymb,latexsym,amsmath,txfonts,color,graphics}
\usepackage[colorlinks,linkcolor=blue,urlcolor=blue,citecolor=black,
plainpages=false,pdfpagelabels,breaklinks]{hyperref}
\usepackage[pdftex]{graphicx}
\usepackage{mathpazo,makeidx}

\usepackage[ paperheight=297mm, paperwidth=210mm,  % or: "paper=a4paper"
             layoutheight=200mm, layoutwidth=120mm,
             layoutvoffset=41.9mm, layouthoffset=45mm, %changed layoutvoffset from 48 to lower text
             centering,
             margin=0pt, includeheadfoot,
             footskip=5mm,
           ]{geometry}

\title{{Quantum Mechanics, Ontology, and Non-Reflexive Logics}\thanks{See the warning note at the end.}}
\author{{{D\'ecio Krause}\thanks{Partially supported by CNPq.}} \\ {Department de Philosophy} \\ {Federal University of Santa Catarina} \\ and \\ Post-Graduate Program in Logic and Metaphysics \\
Federal University of Rio de Janeiro}
\date{{88040-900 Florianópolis, SC $-$ Brazil} \\ \tt{deciokrause@gmail.com}}
\begin{document}
\maketitle

\newtheorem{thm}{Theorem}[section]
\newtheorem{lem}{Lemma}[section]
\newtheorem{cor}{Corollary}[section]
\newtheorem{dfn}{Definition}[section]
\newtheorem{exe}{Example}[section]
\newtheorem{axm}{Axiom}[section]

\newcommand{\ita}{\textit}
\newcommand{\mcal}{\mathcal}
\newcommand{\mbb}{\mathbb}
\newcommand{\mfr}{\mathfrak}
\newcommand{\msf}{\mathsf}
\newcommand{\od}{\mathrm{ord}}
\newcommand{\lra}{\leftrightarrow}
\newcommand{\igual}{:=}
\newcommand{\qst}{$\mathfrak{Q}$}
\newcommand{\Qsim}{\mathrm{Qsim}}
\newcommand{\Sim}{\mathrm{Sim}}
\newcommand{\Proof}{\noindent\textit{Proof:} \,}
\newcommand{\cqd}{{\rule{.70ex}{2ex}}}

\newcommand{\llb}{\llbracket} %with \usepackage{stmaryrd}
\newcommand{\rrb}{\rrbracket}

\thispagestyle{empty}
\pagestyle{empty}

\begin{abstract}
This is a general philosophical paper where I overview some ideas concerning the non-reflexive foundations of quantum mechanics (NRFQM). By NRFQM I mean formalism and an interpretation of QM that considers an involved ontology of non-individuals as explained in the text. Thus, I do not endorse a purely instrumentalist view of QM, but believe that it speaks of something, and then I try to show that one of the plausible views of this `something' is as entities devoid of identity conditions. 

\end{abstract}
\begin{footnotesize}
%\tableofcontents
\end{footnotesize}

\section{Introduction}
A typical way of looking to a scientific theory may be this one:\footnote{Several general characterizations of a scientific theories that apply to most of them were proposed in the literature. In a certain sense, all of them coincide in their main aspects to ours. See for instance \cite{daltor81}, \cite{kraare16}.} it consists of a triple $\mathcal{T} = \langle \mathcal{F}, \mathcal{M}, \mathcal{R}\rangle$, where $\mathcal{F}$ is a mathematical formalism (the mathematical counterpart of $\mathcal{T}$), $\mathcal{M}$ is the class of its models, that is, the `realizations' or the `interpretations' of the theory, and $\mathcal{R}$  represents the set of connection rules which provide the links between the interpretations and the formalism. This is of course just a general scheme that guides us to look to some theories from a very general point of view. Differently from mathematics, physical theories demand the explanation of (at least) one possible realization. Quantum mechanics (QM) would be not different. Below we shall see the details. 

%We can characterize $\mcal{M}$ in two ways. The first one is in accordance with standard model theory: a model of $\mcal{F}$ is a mathematical structure according to which the postulates of $\mcal{T}$ are true. The second one is an informal semantics, described in standard natural language, perhaps supplemented with additional symbols. So, most of 

If we considere a purely instrumentalist point of view (sometimes associated to Bohr), we can say that QM just provides us with ways of computing probabilities. But this is a radical `physical' point of view; in general, philosophers aim at to discuss the kind (or kinds) of world(s) QM tells us about. In doing that, we are faced with interesting and counter-intuitive ways of looking to the objects that form our `quantum reality'. Let me explain this point a little bit, so directing the discussion to that I wish to emphasize. Firstly, I invite you to agree that  Sunny Auyang is right when she says that ``physical theories are about things" (\cite[p. 152]{auy95}) for, if not, we would be speaking of a purely mathematical theory. But, of which kind of things are we speaking about? 

To enlighten the point, let me recall one of Arnold Scharznegger's films, The 6th Day where, coming back to home, he discovers that there is someone  in his place,   a perfect copy of him.\footnote{The story (or history) of  Martin Guerre, a Frenchman who lived in the XVI century and which was also substituted by an imposter, is told by Amanda Gefter in "Quantum mechanics is putting human identity on trial: if our particles have no identity, how can we have?", \href{https://goo.gl/WrC8Dj}{Nautilus}.} The clone act he does, interacts with his wife and sons as he does, puts out the trash as he does, etc. His family does not perceive any difference, and really believe that the person they are interacting with is Adam (the name of the character). But, is he not Adam? To his family, all that imports is that \ita{all happens as if} the person in the house looks as Adam, acts as Adam, etc. That is all. But, you may say, the person in the house \ita{is not} Adam, but just someone quite similar to him. 

Let us go to another example: imagine that someone very rich buys a painting (supposedly) by Picasso and then realize that it is not a legitime Picasso but a copy, a very good one. Okay, you can say: no problem. The copy is not a Picasso, but just a copy. But suppose that Picasso has produced two very similar paintings, one of them acquired by our personage, and that only after the guy has acquired it the second copy appears, and that the specialists are with difficulties to verify the legitimacy of the second painting. What could happen? Probably the the painting (both) would lose their value. Perhaps not, who knows? What imports is that, apparently, identity matters:\footnote{But see below for I will question even this thesis.} independently of the similarities, just one of the guys is Adam, just one of the paintings is the original Picasso and in the case of two Picasso's, something else happens: the value of the painting may decrease. Furthermore, the two paintings can be discerned one each other, for some (yet small) difference in tones, in traces or whatever else would be present. With ordinary objects or our scale, apparently Leibniz's celebrated Principle of the Identity of Indiscernibles holds.

Concerning quantum objects, things are different.  Although we know the no-cloning theorem \cite{woozur82}, all electrons are electrons-Adam, all protons are protons-Adam, and so on. They don't present any intrinsic difference, but are perfectly alike (supposing of course that it makes sense to speak of them as `things' of some kind).\footnote{An useful characterization of quantum particles, adopted by most physicists, was given by E. Wigner in 1939 using group theory \cite{wig59}; see \cite{cas94}.}   Yes, you can say that  certain quantum objects of the same kind, as two electrons (and for other fermions in general), contrarily to bosons, cannot partake the same quantum state since they must obey Pauli's Principle, \ita{hence} they do present differences for according to this principle, no fermions can partake the same quantum numbers in a same situation. But the problem must be put rightly:  think of the two electrons of an Helium atom in its fundamental state. The two electrons are in a superposed state, and yet so they differ by their values of spin in a given direction. But the question is that no one can say which is which.\footnote{And this is not a simple epistemological problem. In assuming this, we should agree in admitting hidden variables of a kind.} Even if we call one of them `Paul' and the other `Peter', before a measurement it is impossible to say which electron is Paul and which one is Peter, contrary to the two supposed Picasso's paintings, for we can write `Peter' in the back of one of them, something we cannot do with electrons.\footnote{A long time ago, Shrödinger stressed that ``you cannot mark an electron, you cannot paint it red'' \cite{sch53}.} Things are worst for bosons, for they may share all their quantum numbers (like in a BEC, a Bose-Einstein Condensate \cite{ket99}), being absolutely indistinguishable by all means provided by the theory and perhaps by any means at all. Identity, here, seems to make no sense. 

But, what identity?  I mean the intuitive idea of a perfect characterization of an object as a sole object. According to our preferred metaphysics (a Leibnizean one), things having identity are always \ita{different} from any \ita{other} things, and can be recognized as such in different situations. Although sometimes we cannot say in what the difference consists of (more on this below), we tend to suppose that it exists. As we shall see soon, this hypothesis is encoded in our metaphysical pantheon. But let me insist on the nature of quantum objects. Think for instance in the methane combustion:
$$\mathrm{CH}_4 + 2 \mathrm{O}_2 \rightarrow \mathrm{CO}_2 + 2 \mathrm{H}_2\mathrm{O}.$$

Here, the four oxygen atoms of the two oxygen molecules will contribute in the reaction so that two of then will form the carbon dioxide molecule and the other two will form the two water molecules. But, which ones? It does not matter.  All oxygen (hydrogen, carbon, electrons, protons, etc.) act the same way in the same circumstances. As another exemple, think of an ionization process of a neutral helium atom. One of its electrons may be realized from the atom in order to form a positive ion. Later, an electron can be captured by the ion in order to form a neutral atom again. What is the difference between the first and the second neutral atoms ou between the realized electron and the absorbed one? None. There are no differences at all! If there were differences, chemistry would not work as it does. A long time ago (1803), John Dalton has put things clear when he stressed that 

\begin{quote}
``[w]hether the ultimate particles of a body, such as water,
are all alike, that is, of the same figure, weight, etc.
is a question of some importance. From what is known,
we have no reason to apprehend a diversity in these particulars:
if it does exist in water, it must equally exist
in the elements constituting water, namely, Hydrogen and
Oxygen. Now it is scarcely possible to conceive how the
aggregates of dissimilar particles should be so uniformly
the same. If some of the particles of water were heavier
than others, if a parcel of the liquid on any occasion were
constituted principally of these heavier particles, it must
be supposed to affect the specific gravity of the mass, a circumstance
not known. Similar observations may be made
on other substances. Therefore we may conclude that the
ultimate particles of all homogeneous bodies are perfectly
alike in weight, figure, etc. In other words, every particle
of water is like every other particle of water, every particle
of Hydrogen is like every other particle of Hydrogen, etc." \cite[pp.142-3]{dal08}
\end{quote}

Quantum objects do not behave as the objects of our scale. They are quite strange. 

\section{Is identity really so fundamental?}
The title of this section is the title of a paper I wrote with my colleague Jonas Arenhart \cite{kraare15}. We discuss and contest and idea that the notion of identity is fundamental, as advanced by O. Bueno \cite{bue14}. I will not repeat the arguments of the paper here, but provide a mix of related ideas instead. 

Standard objects of our surroundings, we believe, \ita{do have identity}, they are \ita{individuals}. This informally means that they have their \ita{own} characteristics which distinguish them from any other object, although a quite similar one, and they can (in principle) be re-identified several times as being \ita{that} individual. This is due to their \ita{identity}. 
This belief is one of the most celebrated metaphysical assumptions of Western philosophy already mentioned above, namely, the Principle of the Identity of Indiscernibles (PII), which has its routs in the antiquity, but was celebrated in Leibniz's philosophy.\footnote{Max Jammer suggests that PII has its roots with the Stoics; see \cite[p.]{jam66}.} PII says that no two objects can have all the same properties, that there are no \ita{solo numero} distinct objects. Objects having the same properties are \ita{indistinguishable}, or \ita{indiscernible}, while identical objects are the very same object. 
PII makes these concepts equivalent. As Rodriguez-Pereyra says, ``the principle states there cannot be \ita{numerically distinct} but perfectly similar things, or that there cannot be \ita{two} perfectly similar things'' \cite[p.15 fn.1]{per14}. Classical logic and standard mathematics, that is, that one which can be build, say, in a standard set theory like the Zermelo-Fraenkel with the Axiom of Choice system (ZFC), incorporates PII in some way; we can say that it is a theorem of classical logic (and classical mathematics). Every object in ZFC (either if it involves or not the \ita{Urelemente}, entities that are not sets but which can be elements of sets --- see \cite{sup72}) is an \ita{individual} in the sense of obeying PII. More precisely, by obeying the standard theory of identity of ZFC, which if formulated as a first-order theory, formalizes the behavior of a binary primitive relational symbol `=' by means of the following postulates: (1) Reflexivity: $\forall x (x=x)$; (2) Substitutivity: $\forall x \forall y (x=y \to (\alpha(x) \to \alpha(y)))$, where $\alpha(x)$ is a formula with $x$ free and $\alpha(y)$ results from $\alpha(x)$ by the substitution of $y$ in some free occurrences of $x$, and (3) the Axiom of Extensionality: $\forall x \forall y (\forall z (z \in x \lra z \in y) \to x=y)$. In order to prove that every object in ZFC is an individual, it suffices to acknowledge that given any object $a$ (represented in ZFC either by a set or by an \ita{Urelement}), we can form the unitary set $\{a\}$ and define the following property, which I call `the identity of $a$', namely, $I_a(x) \lra x \in \{a\}$. It is obvious that the only object obeying $I_a$ is $a$ itself; so, any \ita{other} object will have a difference with $a$ and PII holds. 

Good for mathematics. If in discussing arithmetics `my' number two is different from `yours', we shall have some troubles, yet we both acknowledge that there are different and non equivalent way of defining `two' (for instance, Frege's, Russell's, Zermelo's, von Neumann's definitions, and so on). But in a same context, we can agree that the two twos are the same. But for empirical sciences, I guess that we don't need a so strong assumption; worst, in assuming the standard theory of identity as true for empirical objects in general, mainly in regarding quantum objects, we shall face some problems I will touch on below. Firstly let me recall that long time ago David Hume made things clear (in my opinion) when he discussed the identity of objects in his \ita{Treatise} \cite[Book I, \ita{passim}]{hum85}. In short, he said that there is nothing in an object that justifies our belief that it continues to be itself after two successive observations, the second one after an instant in which the object leaves our field of perception. We can say that it is only a postulate of ours that things happen  this way. Interestingly enough, Schr\"odinger had a similar position when he stressed that 

\begin{quote}
``[w]hen a familiar object reenters our ken, it is usually recognized as a continuation of previous appearances, as being the same thing. The relative permanence of individual pieces of matter is the most momentous feature of both everyday life and scientific experience. If a familiar article, say an earthenware jug, disappears from your room, you are quite sure that somebody must have taken it away. If after a time it reappears, you may doubt whether it really is the same one $-$ breakable objects in such circumstances are often not. You may not be able to decide the issue, but you will have no doubt that the doubtful sameness has an indisputable meaning $-$ that there is an unambiguous answer to your query. So firm is our belief in the continuity of the unobserved parts of the string!" \cite[p.204]{sch98}
\end{quote}

%As in the case of Schwarznegger's film, for all practical purposes (except by a sentimental one, or by an historical fact, say of an ancient Chinese jug) it does not matter either the jug is or not the same. Even in the case of the Chinese jug, or of a painting by Picasso, we can be not absolutely sure that there is not a perfect copy of Guerniza painted by Picasso himself hidden in some place elsewhere. Of course we can doubt that, but it is not logically impossible. Anyway, if instead of earthenware jugs or paintings we are dealing with elementary quantum objects, my claim that it does not matter which particular one is being considered makes perfect sense. To sum up, I think that Hume was in the right direction; after a night in the Museo Reina Sophia, we don't have a logical proof that the Guernica we find in the morning it \ita{identical} to that one we have left  the last evening. Exaggerating a little, the painting it is not the \ita{same}, for it has changed some of its characteristics. But in order to consider this, we would be in need to consider the formation of the notion of an object, something already discussed in the literature (for instance, by Jean Piaget $-$ \cite{pia86}).

It seems that, for physics at least, identity of objects in time cannot be proven logically; for physics at least (not for art or human relationships) what imports is that once we have an object, an indistinguishable one serves as well. So, we need to assume that indistinguishable objects should exist in some way. The problem with this idea is mathematics (and logic) for, as we have seen, PII entails that indistinguishable objects are \ita{the very same}  object; in other words, in ZFC there is not legitimate (\ita{solo numero}) indistinguishable objects. We shall see why next. 

\section{Indistinguishability within classical logic}
Let us consider a structure $\mathcal{E} = \langle S, R_i\rangle$, $i\in I$, built in ZFC.  A structure of this kind may admit automorphisms other than the identity  function (which is of course an automorphism in every structure). For instance, think of a group structure  $\mathcal{G} = \langle G, \star\rangle$. An automorphism of $\mathcal{G}$ is a bijective mapping $h : G \to G$ such that (i) $h(x \star y) = h(x) \star h(y)$, (ii) $h(e) = e$ (where $e$ is the identity element of the group), and (iii) $h(x') = (h(x))'$, where $x'$ is the inverse element of $x$. For instance, take the group defined by $\mcal{Z} = \langle \mbb{Z}, + \rangle$ where $\mbb{Z}$ is the set of the integers and + the standard addition on this set. Then $h : \mbb{Z} \to \mbb{Z}$ defined by $h(x) = -x$ is an automorphism of $\mcal{Z}$, as is easy to prove. If $h$ is an automorphism of the structure $\mathcal{E}$ and $h(a) = b$ for $a, b \in S$, we say that $a$ and $b$ are $\mathcal{E}$-\ita{indistinguishable}. Thus, $2$ and $-2$ are $\mcal{Z}$-indistinguishable. From the point of view of the structure, that is, \ita{from within} the structure, there are no ways of distinguishing between $a$ and $b$: they look the same, for they are invariant by the structures' automorphisms. 
But, are they identical? This of course happens only if there is just one automorphism, namely, the identity function. In this case, we say that the structure is \ita{rigid}.  The interesting thing to recall is that in ZFC every structure can be extended, by adding new relations, to a rigid one. This intuitively means that we can always `go out' the structure and look its elements from the point of view of the outside, and from this point of view, we can realize that the objects we initially thought were indistinguishable, are not indistinguishable at all!  For instance, we can `rigidify' the structure $\mcal{Z}$  by adding the binary relation $<$, the usual linear order on $\mbb{Z}$; that is, the extended structure $\mcal{Z}'= \langle \mbb{Z}, +, < \rangle$ is rigid, as is easy to see. 
The idea of `leaving out' the structure is similar to go to another dimension; in the 1884 book \ita{Flatland} \cite{abo91}, Edwin Abbott created a world in two dimensions, and we can imagine a teenager character of the story who is boring with every one else and decides to keep closed in her room, alone. She supposes nobody can see her. But we, in the third dimension, can. The same happens with our $\mathcal{E}$-indistinguishable objects; looking to the elements of the domain from the extended structure, we can realize that they are not indistinguishable at all. And this can always be done!  \cite{cosrod07} Within the ZFC framework, an object is indistinguishable just from itself and remain so independently of what we `do'\  with them: they are \ita{individuals}. 

The moral of the story is this: within a standard mathematical framework such as ZFC, the only way of considering indiscernible objects is to confine the discussion to a non-rigid structure or something similar which is equivalent to admit as `identical' those objects belonging to a same equivalence class relative to some equivalence relation or to some congruence (other than identity). This is precisely what the formalism of QM does when postulates that only \ita{certain} states are accessible to quantum objects. 
For instance, let us consider two indistinguishable bosons and two possible states $A$ and $B$. The configuration space is the tensor product Hilbert space $\mathcal{H} = \mathcal{H}_1 \otimes \mathcal{H}_2$, where $\mathcal{H}_i$ ($i=1,2$) are the state spaces of bosons 1 and 2 respectively. Note that in order to speak of them, we need to name them, say by calling them `boson 1' and `boson 2'. But these names cannot make them entities with individuality, so we need to provide a mathematical trick in order these names lose their individuation roles. This is done by assuming that vectors like $\psi_1^A \otimes \psi_2^B$ and $\psi_2^A \otimes \psi_1^B$, meaning respectively that particle 1 is at state $A$ and particle 2 is in $B$, and that particle 2 is at state $A$ and particle 1 is in $B$, are not accessible to the particles. Furthermore, these vectors are (in general) different, for the tensor product is not commutative. So, the indiscernibility of the particles cannot be preserved, for it would be a different situation either we consider that particle 1 is at state $A$ and particle 2 is in $B$ or that particle 2 is at state $B$ and particle 1 is in $A$. 
These vectors do not represent possible states for the join system; in the terminology introduced by Michael Redhead, these vectors are examples of \ita{surplus structures}, objects resulting from the formalism that have no physical significance (see \cite{redtel91}). The `right' states are (1) $\psi_1^A \otimes \psi_2^A$, meaning that both are in $A$, (2) $\psi_1^B \otimes \psi_2^B$, meaning that both are in $B$, and (iii) $\frac{1}{\sqrt{2}}(\psi_1^A \otimes \psi_2^B  \pm \psi_2^A \otimes \psi_1^B)$, the plus sign holding for bosons (symmetric states) and the minus sign for fermions (anti-symmetric states --- for fermions, this the only available state, due to Pauli's Principle). The situation (iii) says that \ita{one} particle is in $A$ and that the \ita{another one} is in $B$, but we cannot state which is which. Thus, the formalism preserves Pauli's Principle for fermions; although a permutation of the particles changes the signal of the whole vector, its square remains the same, and that is what imports, for the square of the vector (or wave-function) gives us the probabilities, according to Born's well known rule. 

In doing that, the formalism can be written having ZFC as its mathematical (and logical)  basis. Of course alternative frameworks could be invoked instead, but ZFC suffices and is sufficiently general for the considerations we have in mind. So, we are performing a trick in assuming that only symmetric and anti-symmetric vectors are available for quanta. This is similar to the confinement of the discussion to a non-rigid structure.\footnote{In \cite{domholkra08} --- see also \cite{domholknikra10} --- we have started the development of a formalism that dispenses labels to the particles. This is something to be further explored.}

\section{Ontology} 
The word `ontology' has acquired a number of different meanings in philosophy. Traditional philosophy has qualified it as that part of metaphysics that studies  the general structures of what there is. In this sense, there cannot be distinct ontologies, for what there is is what there is and things,  in the sense of \ita{being}, cannot have two distinct natures.  
But today we have relativized  (or `naturalized') the word to a certain theory or conception; thus we can say that, given a scientific theory, its ontology is described by specifying the kind of entities the theory is compromised with. In this sense, we can say that 
the standard formalism of QM  is compatible with many non-equivalent ontologies  (see \cite{frekra06} for a wide discussion) which can be aggregated to the formalism as a kind of interpretation of it (an intended semantics). 
The standard formalism (Hilbert spaces) does not speak of quantum objects strictly speaking, but just of \ita{states} and \ita{observables}, and we need to provide a parallel discussion (sometimes disliked by physicists) in order to answer simple questions such as `states of what'? 
But, `logically speaking', how can we provide a semantics for the formalism? According to standard semantic procedures (formal semantics), the first step is to define a \ita{domain of discourse}. For instance, we may suppose that we are speaking of a collection of bosons or of another quantum system. Let us fix a collection of indistinguishable bosons (all in the same quantum state). 
Would this domain be a set? Remember what Georg Cantor, the founder father of set theory, said about sets: ``by an `aggregate' (\ita{Menge}) we are to understand any collection into a whole (\ita{Zusammenfassung zu einein Ganzen}) $M$ of \ita{definite and separate objects} $m$ of our intuition or our thought", my emphasis \cite[p.85]{can55}.  In ZFC, due to the Axiom of Extensionality, the set $\{1,1,2,3,3,3\}$ is identical with the set $\{1,2,3\}$ and has cardinal 3.\footnote{Let me recall that there is a \ita{multiset theory} where an element can occur more than once in a multiset, so that $\{1,1,2,3,3,3\}$ nas cardinal 6 \cite{bli88}. But this does not fit QM (see \cite{kra91}), for  the repeated objects are \ita{the very same} object, while in QM no one will agree that in a collection of indistinguishable bosons, they are all the same boson. As we shall see below, words like `the same', `identical', `different' causes troubles here.} So, how to provide an adequate semantics for, say, indistinguishable bosons? 
Hermann Weyl has found a way: he has taken a set $S$, say with $n$ elements, and an equivalence relation $\sim$ defined on $S$. Then the equivalence classes $C_1, \ldots, C_k$ have cardinals $n_1, \ldots, n_k$ so that we have an `ordered decomposition' $n = \sum_{i=1}^k n_i$ (see \cite[p.240]{wey49}) and, as he says, this would be what imports to QM, that is, \ita{the quantity} of elements in each state (equivalence class), and not their individual description. The equivalent classes play the role of characterizing indistinguishable elements of the collection, that is, elements of $S$ that belong to a same equivalence class are taken as indistinguishable. For a discussion of this case, see \cite{kra91}, \cite{frekra06}. But this is a trick, for \ita{we know} that the elements of $S$ are all distinct from one each other by fiat! Leibniz's principle holds, and if we have a finite number of objects, we may even ask for the differences. But, as we know from QM, in certain situations there are none! 

We may say that the differences are in logic. But, is logic measurable? Things became difficult if we try to push deeper this philosophical question. Thus, let me suggest an alternative. Today most logicians and philosophers have no more fears in admitting different kinds of heterodox logics, that is, systems that depart from classical logic in some way. Intuitionistic logic do not accept the general validity of the excluded middle principle; paraconsistent logics do not accept the general validity of the principle of contradiction. But the laws of identity are still taken as a taboo. No one I know accepts to question them, and more, to question the ancient metaphysical rule that there cannot be truly indistinguishable objects. Why not to admit them, at least logically? In the same vein as there is no logical proof that there is no another perfectly similar Guernica hidden somewhere, an ontology composed by truly indistinguishable objects, not made ad hoc by some trick as shown above, looks reasonable. Perhaps such a metaphysics would fit well the claims of QM. So, we are free to try to found a semantics for QM whose domain comprises collections of indistinguishable objects. By the way, as said David Hilbert, the mathematician (and the philosopher) should investigate all possible theories \cite{hil02}. 

\section{The inadequacy of the standard theory of identity}
%Classical mathematics has interesting results. Using the Axiom of Choice, we can prove that the set $\mathbb{R}$ of the real numbers can be well ordered. This means that there is a well ordering such that any non-empty set has a least element according to this order. So, the intervals $[0,1]$ and $[2,3]$ (given in the usual order) have least elements relative to the well ordering. But it is a theorem  that the well ordering cannot be defined by a formula of the language of ZFC; consequently, although we can name the two least elements respectively `Paul' and `Peter', we cannot describe them (by a formula) since we would be in need of the well ordering to do it. Even so, Paul and Peter are distinct real numbers and do present a difference, say Paul $\in [0,1]$ and Peter $\notin [0,1]$. Parenthesis: really, Paul would be described as follows, being $R$ the well ordering: $(\forall x \in [0,1])(R(Paul,x))$, which requires $R$.

Standard theory of identity says that two \ita{distinct} objects present, at least in principle, a difference in their qualities (that is, there is a property obeyed by just one of them). If we have two of them, they are different in this sense, yet we would not be able to specify the distinctive property. For instance, the bosons in a collection of bosons in the same state, if represented in  ZFC, do present differences, yet QM cannot tell us what are they. Remember: in ZFC, if we have a set with more than one element, its elements are \ita{different}. 

Let us fix this typical case, namely, a collection of indistinguishable bosons, say in a Bose-Einstein Condensate (a `BEC').\footnote{A very interesting and didactic page on BECs is \url{http://www.colorado.edu/physics/2000/bec/what_is_it.html}.} A BEC can be obtained by freezing molecules or atoms near to the absolute zero; in such a situation, the objects start acting as if they were just one thing, a `big molecule' \cite{ket99}. They are in the same quantum state and do not present any differences, \ita{but they are not the same object!}. But, let me insist, if we suppose they obey the standard theory of identity, the differences exist, yet not perceptible to us. In this case, we need to assume that there is something (some kind of `variable') hidden in the quantum mechanical description (in the case, a quantum field theoretical description). But no physicist (I suppose) is comfortable in assuming this. The bosons in a BEC do not present any difference, even a hidden one. So, it seems that they should not belong to the class of objects that obey the standard theory of identity. Of course Bohmian quantum mechanics (BQM) says differently, for it agrees with the `classical' metaphysics of individuals. In BQM, all particles have positions that distinguish them one another. The problem is that as Carlo Rovelli says, the particles "do not revel"\ their positions; they are hidden to us (and, we could add, also to the gods) \cite[p.269, fn.55]{rov18}.

\subsection{My proposal}
What, them? In my opinion, the better way to deal with entities of this kind is to separate the two notions that are merged in the standard theory of identity: indistinguishability and identity. Indistinguishable objects share properties; you and me apparently are indistinguishable regarding the `property' \ita{to have interest in quantum physics}. But for sure we have several other differences. `Truly' indistinguishable objects do not present any difference by definition; they share \ita{all} their properties. Identical objects are not distinct objects, but the same one. In other words, there is no more than one object.
In my opinion, the first concept is useful in quantum physics, while the second one causes troubles, and is useful (perhaps) only in mathematics, art and human relationships. In fact, if we assume the standard theory of identity for bosons, we need to assume also the corresponding `theory of difference' for them, which says that if we have two of them, a difference exists. But, which one?
So, out of a purely metaphysical hypothesis, in order to speak of a suitable semantics for quantum languages, we need a mathematical theory were these two concepts are not taken as equivalent. This theory is called \ita{quasi-set theory} we shall see below. But, first, let me say something about the underlying logic of such a theory, a \ita{non-reflexive logic}. 

\section{Non-reflexive logics, and Schrödinger}
Generally speaking, non-reflexive logics are non-classical logics that deviate from classical logic with respect to the notion of identity. Since we can have classical logic \ita{without} identity (see \cite{men79}), in order to characterize them it is necessary to provide a way to modify in some sense the way classical logic  deals with identity. And this may go as follows. If the logic does not contain a primitive binary predicate to be interpreted as identity, we chose a binary predicate of the language and associate to it the diagonal of the domain of the interpretation $D$, namely, the set $\Delta_D = \{\langle x, x \rangle : x \in D\}$. If the logic comprises a primitive binary predicate of identity, we associate the same set to it. The problem is that, in usual parlance, identity (the informal notion of `being the very same') cannot be axiomatized. It means that we can never know if the associated set is in fact the diagonal of the domain or another set characterized as the quotient set of the domain by some congruence relation; the structures are elementary equivalent, what means that the same sentences are true in both of them (for details, see \cite{hod83}, \cite[p.83]{men79}). So, from the point of view of the first order language, we cannot distinguish between the two structures and, then, we never know if the predicate of identity really stands for the diagonal of the domain, that is, either we are speaking of the individuals of the domain or of equivalence classes of them. 

A typical way of departing from the standard notion of identity is to try to violate the Principle of Identity in some way. But there is no \ita{the} Principle of Identity, for it can be formulated in several non equivalent ways. For instance, at the propositional level, we can write $p \to p$ where $p$ is a propositional variable. In this case, we can interpret the implication as \ita{cause}, prescribed by suitable axioms (see \cite{sylcos88}). That is, $p \to q$ is read $p$ \ita{causes} $q$, and it is assumed that nothing can case itself, so, $p \to p$ does not hold. Other systems can be obtained by considering the Principle of Identity as formulated in a first order language, namely, $\forall x (x = x)$, where $x$ is an individual variable, also called the reflexive law of identity. The negation of this rule reads $\exists x (x \not= x)$. But we are not claiming that there is something which is not identical to itself; the principle can be violated simply by assuming that the predicate of identity does not hold in general, that is, that there may exist objects in the domain to which it does  not make sense to say that they are either equal or different. In this case, $x=y$ simply does not have sense. A typical case is Schrödinger's idea that the notion of identity does not have sense for elementary particles in quantum mechanics. This is the case we have in mind. Inspired in this idea, we have developed Schrödinger Logics in which the Principle of Identity in this first-order form does not hold in general (for details and historical references, see \cite[chap.8]{frekra06}). The theory to be presented below incorporates this idea.

Other discussions on non-reflexive logics can be found in \cite{cosbue09}, \cite{kraare18}, \cite{kra94}, and in the references therein.

\section{Quasi-set theory}
 In this section I shall sketch a minimal nucleus of quasi-set theory just to give to the reader a general idea of how it works. Later, I shall say something about `quantum semantics'. 
 In what follows, ZFU stands for the Zermelo-Fraenkel set theory with \ita{Urelemente} \cite{sup72}. Thus, our intended domain comprises different kinds of entities; standard ZFU is compatible with the existence of \ita{sets} and the \ita{Urelemente}, or simply \ita{atoms}. In the theory sketched below, $\mfr{Q}$, there are two kinds of atoms; the $M$-atoms play the role of the \ita{Urelemente} of ZFU (and in the intended semantics stand for the usual objects of our surroundings $-$ at our scale, the \ita{individuals}), while the $m$-atoms have a different behavior, and will be thought of as representing elementary particles (either in non-relativistic quantum mechanics or in quantum field theories; both situations can be covered by the formalism, although we shall be speaking more of quantum mechanics). I will modify a little some already presented versions of the theory, mainly in postulate ($qc_1$), where I have admitted the possibility that a collection may do not have a cardinal (as in the case of quantum field theories, where creation and anhilation operators are introduced).  Another modification is concerning the definition of extensional identity given below (definition \ref{dfn}v).

Let us call $\mathfrak{Q}$ a first order theory whose primitive vocabulary contains, beyond the vocabulary of standard first order logic without identity (propositional connectives, quantifiers, etc. $-$see \cite{men79}), we have the following specific symbols: (1) three unary predicates $m$, $M$, $Z$, (2) two binary predicates $\in$ and $\equiv$, (3) one unary functional symbol $qc$.  Notice once again that identity is not part of the primitive vocabulary, and that the only terms in the language are variables and expressions of the form $qc(x)$, where $x$ is an individual variable, and not a general term.\footnote{This restriction avoids that, for instance, $qc(qc(x))$ turns to be a term.} 
The intuitive meaning of the primitive symbols is given as follows:

%-----------------------------------------------------------------
\renewenvironment{enumerate}{\begin{list}{}{\rm \labelwidth 0mm
\leftmargin 5mm}} {\end{list}}

\begin{enumerate}
\item (i)  $x \equiv y$     ($x$ is indiscernible from $y$)
\item (ii) $m(x)$ \ ($x$ is a `micro-object', or an $m$-atom)
\item (iii) $M(x)$  ($x$ is a `macro-object' or an $M$-atom)
\item (iv) $Z(x)$  ($x$ is a `set' $-$ a copy of a ZFU set)
\item (v) $qc(x)$  (the quasi-cardinal of $x$)
\end{enumerate}

The underlying logic of $\mfr{Q}$ is a kind of \ita{non-reflexive logic}, where the standard theory of identity does not hold (see \cite[chap. 8]{frekra06}, \cite{are14}, \cite{cosbue09}). 
Now, we introduce some definitions, with the intuitive interpretation attributed to them.

\begin{dfn}\label{dfn} \hfill{}
\begin{enumerate}
\item (i) $Q(x) \igual  \neg (m(x) \vee M(x))$ \, ($x$ is a qset)
\item (ii) $P(x) \igual Q(x) \wedge \forall y (y \in x \to m(y)) \wedge \forall y \forall z (y \in x \wedge z \in x \to y \equiv z)$ \\ ($x$ is a pure qset, having only indiscernible $m$-atoms as elements.
\item (iii) $D(x) \igual M(x) \vee Z(x)$ \\ ($x$ is a \textit{Ding}, a `classical object' in the sense of Zermelo's set theory, namely, either a set or a `macro \ita{Urelement}'.)
\item (iv) $E(x)  \igual Q(x) \wedge \forall y (y \in x \to Q(y))$ \\ ($x$ is a qset whose elements are qsets.)
\item (v) $x =_{E} y \igual (Z(x) \wedge Z(y) \wedge \forall z ( z \in x \leftrightarrow z \in y )) \vee (M(x) \wedge M(y)  \wedge \forall_{Z} z (x \in z \leftrightarrow y \in z))$  (Extensional identity)--- we shall write simply $x=y$ instead of $x=_E y$ from now on. Notice that the expression $x = y$, when either $x$ or $y$ is an $m$-atom, yet it can be written, it does not have any meaning in the theory.\footnote{This is similar to name $\mathcal{R}$ the collection $\mathcal{R} = \{x : x \notin x \}$ (Russell's set), which can be expressed in the language of ZFC but is not a \ita{set} of this theory, supposed consistent.} The notion of identity applies just to sets and to $M$-atoms. Furthermore, just to explain the terminology, sometimes I use relativized quantifiers: for instance, $\forall_Q x \gamma$ means $\forall x (Q(x) \to \gamma)$, while $\exists_Q x \gamma$ means $\exists x (Q(x) \wedge \gamma)$; these same for predicates other than $Q$.
\item (vi) $x \subseteq y \igual \forall z (z \in x \to z \in y)$ \, (subqset) \\
Important to realize here the conditional in this definition. Having no identity, we may be in trouble in trying to prove that a certain $m$-atom belongs to a quasi-set, for it should be \ita{identical} to some element of it. This fact does not matter for our purposes. The definition says that \ita{if} $z$ belongs to $x$ \ita{then} $z$ belongs to $y$. In $\mfr{Q}$, it suffices to prove (or to assume) that there is an indiscernible from $z$ in $x$. For instance, in a Litium atom $1s^2 2s^1$, it suffices to say that there is one electron in the outer shell; it does not matter which one.
\end{enumerate}
\end{dfn}

As I have said, $\mathfrak{Q}$ is a theory compatible with the existence two kinds of ur-elements, the $m$-atoms and the $M$-atoms, and also collections formed by either atoms or other collections, the qsets, or by both, atoms and qsets. The theory does not postulate the existence of atoms, as in the standard presentations of ZFU. Some qsets are specially important: when their transitive closure does not contain $m$-atoms,  they contain only what we call `classical objects' of the theory (objects satisfying $D$); items fulfilling this condition satisfy the predicate $Z$ and  coincide with the sets in ZFU. So, classical mathematics can be built inside $\mathfrak{Q}$, in its \ita{classical part}.

The main idea motivating the development of the theory is that some items are non-individuals (roughly speaking, entities for which the standard notion of identity does not apply), and does not obey the notion  encapsulated in the definition of extensional identity. As explained above this concept is not defined for $m$-atoms, the items which intuitively represent quantum indistinguishable objects. Thus, on one side, these things `do not have identity', that is, it does not make sense to say they are identical or different and, on the other side, the indistinguishability relation holds for every item of the theory, so $m$-atoms may be indistinguishable without being identical. Important to notice that in saying that some entities are non-individuals, we are not supposing that we cannot speak of them; really, we \ita{can} speak of them, that is, we can \ita{write} $x=y$ even for $m$-atoms, but this expression does not have a sense according to the theory: it says nothing. For instance, a qset of indiscernible $m$-atoms may have a quasi-cardinal greater than one, say 5, and so we can think of five entities in some situation, although they cannot be discerned then in any way. Furthermore, quantified expressions must be interpreted adequately; again in considering a qset of indistinguishable objects (say, bosons in a BEC), the universal quantifier says `all elements of the BEC', while the existential quantifier says `some element of the BEC'. Thus, universal quantification does not mean `each' element of the qset (which would presuppose identity) as the standard interpretation suggest (for more on this point, see \cite{kraare15}).

\subsection{The postulates of ${\mathfrak Q}$}
Besides postulates for classical first-order  logic without identity (which we shall not list here), we introduce the specific postulates for $\mathfrak{Q}$.

\begin{enumerate}\label{postulatesindist}
\item ($\equiv_1$) \, $\forall x (x \equiv x)$
\item ($\equiv_2$) \, $\forall x \forall y (x \equiv y \to y \equiv x)$
\item ($\equiv_3$) \, $\forall x \forall y \forall z (x \equiv y \wedge y \equiv z \to x \equiv z)$
\item ($=_4$) \, $\forall x \forall y (x = y \to (\alpha(x) \to \alpha(y)))$, with the usual restrictions.
\end{enumerate}

The first three postulates say that indistinguishability is an equivalence relation. Now, this relation is not necessarily compatible with  other primitive predicates; so we can keep identity and indistinguishability as distinct concepts. In fact, if $x$ and $y$ are indistinguishable $m$-atoms and being  $z$ a qset,  $x \in z$ does not entail that $y \in z$, and conversely. The fourth postulate says that substitutivity holds only for identical things, that is, for `classical' things. 

\paragraph{Remark:} Someone may say that we are presupposing identity in the metalanguage when we say that variables $x$ and $y$ are different. This is true, but does not collapse the theory. We have a similar situation for instance in paraconsistent logics \cite{coskrabue06}, which are logics apt to deal with contradictory sentences. That is, the Principle of Contradiction in the form $\neg (\alpha \wedge \neg\alpha)$ does not hold in general. But, in elaborating such systems, we do use the principle as being true, for no one would suggest that something is a formula \ita{and} is not a formula, say. In other words, nothing is a formula and not a formula at once. This `constructive' character of scientific theories (and of logics) is discussed in the Chapter 3 of \cite{kraare16}. 

Other postulates are:
\begin{enumerate}
\item ($\in_1$) $\forall x \forall y (x \in y \to Q(y))$  \\ If something has an element, then it is a qset; in other words, the atoms have no elements (in terms of the membership relation).
\medskip
\item ($\in_2$) $\forall_D x \forall_D y (x \equiv y \to x = y)$ \\ Indistinguishable \textit{Dinge} are extensionally identical. This makes $=$ and $\equiv$ coincide for this kind of entities.
\medskip
\item ($\in_3$) $\forall x \forall y [(m(x) \wedge x \equiv y \to m(y)) \wedge (M(x) \wedge x=y \to M(y)) \wedge (Z(x) \wedge x=y \to Z(y))]$
\medskip
\item ($\in_4$) $\exists x \forall y (\neg y \in x)$ \\ This qset can be proved to be a set (in the sense of obeying the predicate $Z$), and it is unique, as it follows from the axiom of weak extensionality we shall see below. Thus, from now own we shall denote it, as usual, by `$\emptyset$'.
\medskip
\item ($\in_5$) $\forall_Q x (\forall y (y \in x \to D(y)) \lra Z(x))$ \\
This postulate grants that something is a set (obeys $Z$) iff its  transitive
closure does not contain $m$-atoms. That is, \textit{sets} in \qst\ are those
entities obtained in the `classical'  part of the theory.
\medskip
\item ($\in_6$) $\forall x \forall y \exists_Q z (x \in z \wedge y \in z)$  
%This axiom says that, given $x$ and $y$, there exists a qset $z$ having them as elements. The problem is that there may be other elements in $z$ as well. But, with the Separation Schema below, we can `separate' from $z$ a qset using the formula $\alpha(w) \lra w \equiv x \vee w \equiv y$. So, this qset has one indistinguishable from $x$ or $y$, or from both. This is the `pair' but, let me remark, it may have more than two elements (that is, its quasi-cardinal may be grater than 2). More on this below.
\end{enumerate}

\begin{enumerate}
\item ($\in_7$) If $\alpha(x)$ is a formula in which $x$ appears free, then $$\forall_Q z \exists_Q y \forall x (x \in y \lra x \in z \wedge \alpha(x)).$$
This is the  Separation Schema. We represent the qset $y$ as follows:   $$[ x \in z : \alpha(x)].$$ When this qset is a set, we write, as usual, $\{ x \in z : \alpha(x) \}.$
\medskip
\item  ($\in_8$) $\forall_Q x (E(x) \to \exists_Q y (\forall z (z \in y \lra \exists w (z \in w \wedge w \in x)))$. \\ The union of $x$, written $\bigcup x$. Usual notation is used in particular cases.
\end{enumerate}

\subsection{Some basic concepts}

From ($\in_6$), by the Separation Scheme using $\alpha(w) \lra w \equiv x \vee w \equiv y$, we get a subqset of $z$ which we denote $$[x,y]_z$$ which is the qset of the indiscernibles of either $x$ or $y$ that belong to $z$. When $x \equiv y$, this qset reduces to $$[x]_z$$ called the qset of the indiscernibles from $x$ that belong to $z$. The qset $[x,y]_z$ does not have necessarily only \ita{two} elements (that is, we may have $qc([x,y]_z) > 2$), for there may be more than just one indistinguishable from $x$ or $y$ in $z$. Given the qset $z$ and one of its elements, $x$, the collections $[x]$ and $[x]_z$ stand for \textit{all} indiscernible from $x$ and the qset of the indiscernible from $x$ that belong to $z$ respectively. (Usually, $[x]$ is too big to be a qset $-$ as in general are collections of \ita{all} objects so and so, as in standard set theory.)

Later, with the postulates of quasi-cardinals, we will be able to prove $[x]_z$ has a subqset whose quasi-cardinal equals to 1, written
$$\llb x \rrb_z.$$

We call it the \ita{strong singleton} of $x$ (really, \textit{a} strong singleton of $x$, for we cannot grant that it is unique). It has just one element, and we can think of this element \textit{as if} it were $x$, but  it follows from the definition that all we can know  is that $\llb x \rrb_z$ contains \emph{one object of the `species'} $x$. That is, $qc(\llb x \rrb_z)=1$, so there is one item indistinguishable from $x$ in this qset. To prove that this element is $x$, we need identity.

\subsection{Other postulates and definitions}

\begin{enumerate}
\item ($\in_9$) $\forall_Q x \exists_Q y \forall z (z \in y \lra w \subseteq x)$, \\ The power qset of $x$, denoted $\mathcal{P}(x)$.
Interesting here is that we would be in trouble to teach quasi-set theory to children. For instance, take a qset $x$ with cardinal 2 so that its elements (call them $y$ and $z$) are indistinguishable. Now try
to write the qset $\mathcal{P}(x)$. You can't do it significantly. Really, it results that the two subsets with quasi-cardinal 1 are indistinguishable (by the Weak Extensionality Axiom), so, something like $\mathcal{P}(x) = [\emptyset, [y], [z], x]$ has no clear sense. Even so, as we shall see from axiom ($qc_7$), the quasi-cardinal of $\mathcal{P}(x)$ is 4.

\medskip
\item ($\in_{10}$) $\forall_Q x (\emptyset \in x \wedge \forall y (y \in x \to y \cup [y]_x \in x))$, \\ The infinity axiom.
\medskip
\item ($\in_{11}$) $\forall_Q x (E(x) \wedge x \not= \emptyset \to \exists_Q y (y \in x \wedge y \cap x = \emptyset))$, \\  The axiom of foundation, where $x \cap y$ is defined as usual.
\end{enumerate}

\begin{dfn}[Weak ordered pair]
\begin{equation} \langle x, y \rangle_z \igual [[x]_z, [x,y]_z]_z
\end{equation} \end{dfn}

Then, $\langle x, y \rangle_z$ takes all indiscernible from either $x$ or $y$ that belong to $z$, and it is called the `weak' ordered pair, for it may have more than two elements. Sometimes the sub-indice $z$ will be left implicit.

\begin{dfn}[Cartesian Product] Let $z$ and $w$ be two qsets. We define the cartesian product $z \times w$ as follows:
\begin{equation}
z \times w \igual [\langle x, y \rangle_{z \cup w} : x \in z \wedge y \in w] \end{equation} \end{dfn}

Functions and relations cannot also be defined as usual, for when there are $m$-atoms involved, a mapping may not distinguish between arguments and values. Thus we provide a wider definition for both concepts, which reduce to the standard ones when restricted to classical entities. Thus,

\begin{dfn}[Quasi-relation] A qset $R$ is a binary quasi-relation between to qsets $z$ and $w$ if its elements are weak ordered pairs of the form $\langle x , y \rangle_{z \cup w}$, with $x \in z$ and $y \in w$. \end{dfn}

\begin{dfn}
[Quasi-function] f is a quasi-function among q-sets $A$ and $B$ if and only if f is quasi-relation between $A$ and $B$ such that for every $u \in A$ there is a $v \in B$ such that if $\langle u, v\rangle \in f$ and $\langle w,z\rangle \in f$ and $u \equiv w$ then $v \equiv z$. \end{dfn}

In words, a quasi-function maps indistinguishable elements to indistinguishable elements. An interesting question concerns the more specific kinds of functions, that is, injections, surjections and bijections. One can, with some restrictions, define the corresponding concepts, but we shall not present them here (see \cite[chap. 7]{frekra06}).

%-----------------------------------------------------------
\subsection{Postulates for quasi-cardinals}
Notice that in \qst\ the standard notion of identity is not defined for some entities (definition \ref{dfn}v). Now, the identity concept is essential to define many of the usual set theoretic concepts of standard mathematics, such as well order, the ordinal attributed to a well ordered set, and the cardinal of a collection. Since  identity is to be senseless for some items in \qst, how can we employ these notions? One alternative would be to look for different formulations employing methods that do not rely on identity. Another possibility would be to introduce these concepts as primitive and give adequate postulates for them. Concerning the notion of cardinal, there are interesting issues we should acknowledge. First of all, in \qst, there cannot be  well-orders on quasi-sets of  indistinguishable $m$-atoms. Really, a well-order  would imply, for example, that there is a least element relative to this well order, a notion which could only be formulated if identity was defined for $m$-atoms, for this element would be different from any \ita{other} element in the quasi-set. Second, the usual claim that aggregates of quantum entities can have a cardinal but not an ordinal demands a distinction between the notions of ordinal and of cardinal  of a quasi-set; this distinction is made in \qst\ by the introduction of cardinals as a primitive notion, called quasi-cardinals.\footnote{As shown by Domenech and Holik, we can define quasi-cardinals for finite qsets in \qst, without resulting that the qset will have an associated ordinal in the usual sense; see \cite{dom07}.}

Let us see the postulates for quasi-cardinals; for details and motivations, see \cite[Chap.7]{frekra06}, \cite{frekra10}. Here $\alpha$,  $\beta$, $\ldots$ stand for cardinals (defined as usual in the classical part of the theory, that is, in the theory \qst\ when we rule out the $m$-atoms):

\begin{enumerate}
\item ($qc_1$)\label{axi} $\forall_Q x (\exists_Zy (y=qc(x)) \to \exists ! y (Cd(y) \wedge y=qc(x) \wedge (Z(x) \to y= card(x)))$ \\  In words, if the qset $x$ has a quasi-cardinal, then its (unique) quasi-cardinal is a cardinal (defined in the `classical' part of the theory) and coincides with the cardinal of $x$ stricto sensu if $x$ is a set. As recalled above, this axiom does not grant that every qset has a well defined quasi-cardinal.

\medskip
\item ($qc_2$) $\forall_{Q} x (\exists y(y = qc(x) \to x \not= \emptyset \to
qc(x) \not= 0)).$ \\ Every non-empty qset that has a quasi-cardinal has a non-null quasi-cardinal.

\medskip
\item($qc_3$) $\forall_{Q} x (\exists_Z \alpha (\alpha = qc(x)) \to  \forall \beta (\beta \leq \alpha \to \exists_{Q} z (z \subseteq x \wedge qc(z) = \beta)))$ \\  If $x$ has quasi-cardinal $\alpha$, then for any cardinal $\beta \leq \alpha$, there is a subqset of $x$ with that quasi-cardinal.
\end{enumerate}

In the remaining axioms, for simplicity, we shall write $\forall_{Q_{qc}} x$ (or $\exists_{Q_{qc}} x$) for quantifications over qsets $x$ having a quasi-cardinal.

\medskip
\begin{enumerate}
\item ($qc_4$) $\forall_{Q_{qc}} x \forall_{Q_{qc}} y (y
\subseteq x \rightarrow qc(y) \leq qc(x))$

\medskip
\item ($qc_5$) $\forall_{Q_{qc}} x \forall_{Q_{qc}} y
(Fin(x) \wedge x \subset y \to qc(x) < qc(y))$
\end{enumerate}

It can be proven that if both $x$ and $y$ have a quasi-cardinal, then $x \cup y$ has a quasi-cardinal. Then,

\medskip
\begin{enumerate}
\item($qc_6$) $\forall_{Q_{qc}} x \forall_{Q_{qc}} y (\forall w (w \notin x \vee w \notin y) \to qc(x \cup y) = qc(x) + qc(y))$
\end{enumerate}

In the next axiom,  $2^{qc(x)}$ denotes (intuitively) the quantity of subquasi-sets of $x$. Then,

\medskip
\begin{enumerate}\label{2qc}
\item ($qc_7$)   $\forall_{Q_{qc}} x (qc(\mathcal{P}(x)) = 2^{qc(x)})$
\end{enumerate}

This last axiom enables us to think of subqsets of a given qset in the usual sense; for instance, if $qc(x) = 3$, the axiom says that there exists $2^3 = 8$ subqsets, and axiom $(qc_3)$ enables us to think that there are subqsets with 0, 1, 2 and 3 elements. Furthermore, as we have seen above, in \qst\ we can prove that given any object $a \in z$ (either an $m$-atom, $M$-atom or quasi-set) we may obtain the strong singleton of $a$, $\llb a \rrb_z$ whose quasi-cardinal is 1. Important to insist that there is no sense of saying, within \qst, that $a$ is the only element of $\llb a \rrb$, for in order to prove it we need identity. Anyway, \qst\ is consistent with this idea and we may reason \ita{as if} this is really so. So, we can think that within \qst\ that we may have \ita{a certain $m$-atom}, without identifying it, except that it has some characteristics or properties, and not others (for instance, it may be discernible from another $m$-atom $b$). That $m$-atoms may have different properties can be seen from the fact of \qst\ that \qst\ \ita{doesn't prove} the Substitutivity of Indiscernibles, that is,
$$ \mathfrak{Q} \not\vdash a \equiv b \to \forall_Q z (a \in z \lra b \in z).$$

To prove this result, it suffices to take $\llb a \rrb_z$. Since $qc(\llb a \rrb_z) = 1$, $a$ and $b$ cannot belong both to this qset, except if $a$ is identical to $b$, which cannot be assumed in the case of $m$-atoms. So, in an extensional context (and \qst\ is also a kind of an extensional theory, although this should be qualified), we can read $a \in z$ as $a$ having a certain `property' (whose `extension' would be $z$). So,  even indistinguishable $m$-atoms may have distinct properties, as the two electrons in an Helium atom in its fundamental state have different values of spin in a given direction. As for bosons, let $a$ and $b$ name two bosons in a BEC. The first think to acknowledge is that these names do not make sense, for they cannot individualize the named bosons; furthermore, in \qst, the strong singletons $\llb a \rrb_z$ and $\llb b \rrb_z$ are indistinguishable (by the Weak Extensionality Axiom below), hence there are no differences among them (yet they are not \ita{the same} quasi-set). 

%==================================
\subsection{The  Weak Extensionality Axiom}
The weak extensionality axiom generalizes the usual extensionality axiom. Intuitively, it grants us that two q-sets with the same quantity of the same kinds of elements are indistinguishable. For that, we need two extra definitions, the notion of similarity between q-sets, denoted by $\Sim$, and the notion of Q-similarity, denoted $\Qsim$. Intuitively speaking, similar q-sets have elements of the same kind, and q-similar q-sets have elements of the same kind, and in the same quantity:

\begin{dfn} \hfill{}
\begin{enumerate}
\item (i) $\Sim(x, y) \igual \forall z \forall w (z\in x \wedge w \in y \rightarrow z \equiv y)$;
\item (ii) $\Qsim(x, y) \igual \Sim(x,y) \wedge qc(x) = qc(y)$.
\end{enumerate}
\end{dfn}

The weak extensionality axiom reads as follows:

 \medskip
\noindent ($\equiv_{12}$)  $\forall_{Q} x \forall_{Q} y ((\forall z (z \in x/{\equiv} \to
\exists t (t \in y/{\equiv} \wedge  \Qsim(z,t))))
\wedge \forall t (t \in y/{\equiv} \to \exists z
(z \in x/{\equiv} \wedge \wedge \Qsim(t,z))) \to x \equiv y)$

\medskip
Intuitively speaking,  qsets that have `the same quantity' (given by their q-cardinals) of elements of the same kind are indiscernible. 

The following theorem express the invariance by permutations in \qst, and with this result we finish our revision. To prove it, we shall assume another result, namely, that $y \subseteq t$ entails $qc(x - y) = qc(x) - qc(y)$; let us call this result Theorem ($\star$) (the proof can be found in \cite[chap.7]{frekra06}). The theorem goes as follows: 

\begin{thm}[Invariance by Permutations] Let $x$ be a finite qset such that
$\neg(x = [z]_t)$ for some $t$ and let $z$ be an $m$-atom such that $z \in x$. If $w \in t$, $w \equiv z$ and $w \notin x$, then there exists $\llb w \rrb_t$ such that
$$(x - \llb z \rrb_t) \cup \llb w \rrb_t \equiv x$$
\end{thm}\label{unobser}%\label
\Proof  Case 1: the only element of $\llb z \rrb_t$ does not belong to $x$. Then $x - \llb z \rrb_t = x$. Let  $w$ be so that its only element belongs to $x$ (for instance, it may be $z$). Then $(x - \llb z \rrb_t) \cup \llb w \rrb_t = x$, hence the theorem. Case 2: the only element of $\llb z \rrb_t$ belongs to $x$. Then $qc(x - \llb z \rrb_t) = qc(x) -1$ by the mentioned theorem ($\star$). Let $\llb w \rrb_t$ be such that its only element is $w$ itself, so $(x - \llb z \rrb_t) \cup \llb w \rrb_t = \emptyset$. Hence, by Postulate ($qc_7$), $qc(x - \llb z \rrb_t) = qc(x)$. Thus, by the Weak Extensionality Axiom, the theorem follows.
For more details, see \cite[chap.7]{frekra06}. \cqd

\bigskip
 In words, two indiscernible elements $z$ and $w$, with $z \in x$ and $w \notin x$, expressed by their strong-singletons $\llb z \rrb_t$ and $\llb w \rrb_t$, are `permuted' and the resulting qset $x$ remains indiscernible from the original one. The hypothesis that $\neg(x = [z]_t)$  grants that there are indiscernible from $z$ in $t$ which do not belong to $x$. This theorem has a `physical' interpretation: the qset $x$ must be a neutral atom which is to be ionized by realizing an electron in order to become a negative ion. Thus the $m$-object $z$ would represent an electron in the outer shell, while $w$ is `another' electron not in the atom (these words are to be understood metaphorically). Thus, the electron $z$ is realized and, in another experiment, an electron is captured again so that the atom becomes neutral again. The question is: is this last neutral atom \ita{the same} (identical) to the first one? Of course, this would be so if and only if the captured electron is, \ita{ceteris paribus}, exactly the same as the realized one.  But,  is there any sense in saying that the realized electron is \ita{identical} with the captured one? Quasi-set theory escapes from this dilemma by assuming that the basic notion is that of indiscernibility; the electrons are indiscernible, so as the neutral atoms. And this is enough for physics. Philosophically, we have  again our main thesis: the notion of identity is just a useful notion, not an essential one (in \cite{kraare15}, we discuss this thesis with some care). 

\section{A sketch of a quasi-set semantics for quantum languages}
The word `semantics' has also acquired a lot of meanings (as `ontology' did, as we have seen). In the context that interests us here, it refers to the possible links between a certain language, in general a formal language, and certain `realities' that lie outside the language. In other words, the task is to attribute \ita{meaning} to certain terms of the language in order to say that by using its resources we can \ita{speak} of certain entities that form the \ita{domains of application} of the language. Of course a formal language, and the logic of axiomatized theories in general, can make reference to infinitely many different domains. Usually we chose one of them to be our \ita{intended interpretation}. For first order languages we have a rather well developed theory involving semantics in this sense, Model Theory. But for more sophisticated languages (higher-order languages) there is no a general theory of its \ita{models}, that is, interpretations where the postulates of the theory are \ita{true}. By `more sophisticated languages' I mean mainly the languages of physical theories, whose models are in general not first-order structures or, as I prefer to call them, \ita{order-1} structures, composed by a domain (or several domains) and relations and operations relating and operating with the elements of these domains only. In physics (and even in mathematics), in general the relations and operations we have relate (and operate) also with sets of such elements, with functions and matrices formed with them, and so on. These \ita{order-n} structures ($n >1$) need to be dealt with case by case. This is particularly so with QM. So, let us analyse minimally how we can deal with this question in the scope of quasi-set theory.\footnote{Updated: in 2016, Jonas Arenhart and I published a book where more details on this discussion are given; see \cite{kraare16}.} 

In this section I shall sketch a minimal semantics for QM build in $\mfr{Q}$. I will postpone to another opportunity the details and the formal proofs. I will simply justify the general idea of using quasi-sets to make sense the existence of indistinguishable but not identical objects. In speaking of `quantum languages', we need to take some care. Yuri Manin has recalled that quantum mechanics has no its \ita{own} language, making use of a fragment of the language of standard functional analysis (the Hilbert space formalism) \cite[p.80]{man10}. But, inspired by Dalla Chiara and Toraldo di Francia \cite{daltor93} and by G. Cattaneo \cite{cat93}, we can suppose a suitable language  incorporating the standard logical vocabulary, plus the following nonlogical symbols:
\begin{enumerate}
\item (i) A collection of monadic predicates $P_i$ ($1=1, \ldots, n$) to represent `meaningful properties' of the quantum systems, the \ita{observables}. A typical case may be `the value of the spin in the $z$-direction'. Let me call $\msf{P}$ this collection.
\item  (ii) The quantum systems are referred to in different instants of time (according to the intended interpretation) by individual parameters (generic names) $a_1, a_2, \ldots, a_m$. A parameter acts here as when we write the equation of a straight line in Analytic Geometry as $ax + by +c = 0$ and say that $a, b, c$ are parameters ranging on real numbers; we are not specifying particular numbers, but just making reference to them. So, the parameters of our language simply refer to quantum systems without naming them; these parameters are other kind of variables. It seems clear that the more interesting situations are when the quantum systems for a collection of indiscernible entities, like a BEC. In cases like this one is that $\mfr{Q}$ is useful.
\item (iii) A ternary functional symbol $\mcal{P}$ to be interpreted as \ita{probability} in a way to be described below. 
\end{enumerate}

The semantics goes as follows. In the classical part of $\mfr{Q}$, we can consider all sets referred to below, as for instance the set $\mcal{B}(\mbb{R})$ of the Borelians of the real number line. Thus we consider the following structure as our \ita{quantum structure}:

$$\mathcal{QM} = \langle S, \{H_i\}, \{A_j\}, \{U_k\}, \mcal{B}(\mbb{R})\rangle,$$

\noindent where $S$ is a quasi-set suitable for representing the  quantum systems, the $H$s are Hilbert spaces, the $A$s are Hermitian operators defined on suitable $H$s, the $U$s are unitary operators which provide the dynamic of the system (Schrödinger's equation), and $\mcal{B}(\mbb{R})$ is the set of all Borel  sets  of the real number line. The rules of interpretation are defined as follows:

\begin{enumerate}
\item (i) The parameters $a_1, a_2$  etc. are interpreted either as elements of $S$ or as subsets of $S$. Let me observe that in assuming the structure above, we are introducing explicitly the quantum systems in the semantic considerations, something that is omitted in the usual approaches (see \cite{daltor93}, \cite{cat93}, \cite{vf75}, \cite[pp.203ff]{ish95}). 
\item (ii) Each element $s \in S$ (or $s \subseteq S$) is associated to a unitary Hilbert space $\mcal{H} \in \{H_i\}$, represented by a unitary vector $|\psi\rangle \in \mcal{H}$. Composed quantum systems (when $s \subseteq S$) are associated to tensor product of Hilbert spaces, as usual.
\item (iii) The predicates $P_j$ are associated to Hermitian operators $\hat{A} \in \{A_j\}$ of the attributed Hilbert space. The eingenvalues of $\hat{A}$ are the possible outcomes of the measurements made on $P_j$.
\item (iv) For any triple $\langle s, P, \Delta \rangle \in S \times \msf{P} \times \mcal{B}(\mbb{R})$, we have that $\mcal{P}(s, P, \Delta) \in [0,1]$, and this number is interpreted as the probability that for the physical system in state $|\psi\rangle$ (associated to $s$), a measurement made in $P$ gives a value in $\Delta$. As in the standard formulations of QM, we can assume the notion of probability as given in some suitable way (yet the topic is controversial, as is well known). The postulates describing the behavior of $\mcal{P})$ are those of \cite[pp.62ff]{mac63}. 
\end{enumerate}

Then we can proceed as usual with the quantum formalism and intended interpretation, but now with the quantum systems playing a formal role in the developments. Of course more should be said, but I think that the general idea of using quasi-set theory is done. Anyway, some simple question can be envisaged, as the following ones.

%\begin{enumerate}
%\item (a) To each quantum system $s \in S$, we associate a unitary Hilbert space $\mcal{H} \in \{H_i\}$. Composed quantum systems are associated to tensor product of Hilbert spaces, as usual.
%\item (b) Each observable (physical magnitude) $A$ is associated to a Hermitian operator $\hat{A} \in \{A_j\}$ of the attributed Hilbert space. The eingenvalues of $\hat{A}$ are the possible outcomes of the measurements made on $A$. 
%\item (c) If $\{|\alpha_n\rangle\}$ is an orthonormal basis for $\mcal{H}$ formed by eingenvectors of $\hat{A}$, and if $|\psi\rangle$ is a unitary vector of $\mcal{H}$, then we can write $|\psi\rangle = \sum_n x_n |\alpha_n\rangle$, with $x_n = \langle\alpha_n | \psi\rangle$. If $a_1, a_2, \ldots$ are the eingenvalues associated to the vectors of the basis, then the probability that the measurement of the observable $A$ for $s$ in the state $|\psi\rangle$ be $a_m$ is $|x_m|^2$. I am considering here just non-degenerate operators. 
%\item (d) After a measurement, the state of the system collapses to one of the vectors of the basis, say $|\alpha_t\rangle$ with probability $|x_t|^2$. 
%\item (e) If at time $t_0$ the system was in state $|\psi(t_0)\rangle$, at time $t$ it will be at state $|\psi(t) = U(t) |\psi(t_0)\rangle$; this is the time dependent Schrödinger equation. 
%\end{enumerate}

\subsection{Questions}
Some questions are in order, and then we shall see why I have proposed the use of quasi-sets in the semantics presented above. As Dalla Chiara and Toraldo  di Francia have emphasized, there are two basic questions to be answers by any semantics for such a language, namely (here adapted):

\begin{enumerate}
\item (i) Does $S$ or some of its sub collections  determine a set of $m$ elements in the standard set-theoretical sense?
\item (ii) Can a parameter $a_i$ determine a well defined element of $S$?
\end{enumerate}

As they conclude, ``\ldots both these questions have a negative answer" \cite[p.276]{daltor93}. The reasons are easy to find. If indistinguishable, a collection of quantum systems should not be taken as a set of a standard set theory (recall that any attempt in this sense will be in need of some mathematical trick); secondly, the denotation function (the rules of interpretation) cannot be defined as a standard function, for it will not distinguish the elements of $S$ to which attribute the $a_i$. More could be said, but this is enough to sustain our argumentation that quasi-set theory provides a more suitable framework to developed a semantics for quantum languages.

As I have said, all of this of course deserve further explanations and details. But it is clear that QM works as usual, but now with a more suitable form of semantics which serves to make sense of the informal claims about quantum systems. Note that $S$ may have no a well-defined cardinal number, so this semantics may fit also the case of relativistic quantum mechanics, at least for free fields, that is, the Hilbert spaces can be taken as Fock spaces (see \cite[chap.9]{frekra06}). In \cite{domholkra08} and in \cite{domholknikra10}, more is said about the development of QM in $\mfr{Q}$.

\section*{Warning note}
This is a revised version of a paper with the same name that was written by invitation to be published in a book titled \ita{The Mammoth Book on Quantum Mechanics Interpretations}, edited by Open Academic Press, Berlin, and having as editor a certain Ulf Edvinsson, who has invited me. The book was announced in the page of OAP and should appear by 2016. This never happened. Later I discovered that OAP is in a list of predatory editorial houses and that "Ulf Edvinsson" is (apparently) a fake name. Furthermore, I couldn't contact anyone responding by OAP to retire my name from the announcement of the book and for impeding them to publish the paper. I strongly apologize for such a fault, which is completely mine. Since the subject presented here has been among my preoccupations ever since I met Chico Doria for the first time (in 1987), it is a pleasure to dedicate the stuff to him. And of course I thank  the editors for accepting this version of the paper for this book.

\end{document}